# Multirate control strategies for avoiding sample losses. Application to UGV path tracking.

J. Salt, J. Alcaina, A. Cuenca, A. Baños

**Abstract**

When in a digital control strategy there are samples lost due to limitations, diff t multirate (MR) control options can be used for solving the problem: Dual-rate inferential control (IC) and model-based dual-rate control (MBDR). The objective of this contribution is to analyze, compare, and to assess their behavior under diff t perspectives. Is a dual-rate inferential control better than a model-based dual-rate control? Both options lead to a periodically time-varying discrete-time system and for this reason a lifted modeling is considered. An efficient algorithm is used for computing a MR system's frequency response for these control structures. The robust performance and disturbance effects are studied in detail under sample losses and process uncertainty, and some considerations are reported. A new QFT (quantitative feedback theory) procedure for dual-rate systems analysis is also described. Analysis and simulation examples and experimental results for UGV path tracking are introduced in this work, revealing that MBDR outperforms IC when the model contains important uncertainties.

1. **Motivation**

In many computer control applications it is not possible to maintain an ideal sampling period for the measurement of the variable to be controlled. Various economic or technical conditions can cause instants of blank data to occur. It is a scenario that usually appears in industrial plants but it is also given in mobile robots or precision motion control operations. This loss of samples with respect to the ideal measurement frequency leads to what is called slowly sampled systems. There are several alternatives for the analysis and proposal of solutions in this type of systems. One way is what is called soft sensors [1], with the objective of supplying the infrequent measure using secondary variables (process quality variables or key indicators) that can be used in the control and that usually entail an estimation of the primary variable [2, 3]. When the feedback uses estimated outputs, the control method may also be called "inferential control". Soft sensors can include information processing procedures so as to reduce the measurement time and to eliminate delays, to avoid drift and noise, to detect malfunctioning, to minimize maintenance costs, and some other benefits that make them a key idea in modern control systems [4]. From some control structures introduced many years ago, where the use of both primary and secondary variables on nested and cascaded parallel control loops were used [5, 6], to more recent applications such as [2, 7, 8, 9], all of them have shown relevant industrial practical results. The soft sensor development can be considered using an input-output [10], or state-space model [11]. This last contribution provides the skeleton of the so-called "generalized inferential control" for multirate systems, separating the design of an estimator and that of a compensator in order to study diff t $H_2$ optimal controllers in a state-space framework. Within the soft sensors fi diff t methods such as model-based approaches using black-box modeling, inferential-based proposals, and identification-based solutions have been used. Identification-based methods consider statistical methods such as partial least squares, principal component analysis, independent component analysis, support vector machines and Bayesian methods.

It is worth to point out that there are two tendencies with respect to the problem addressed in this contribution: one based on measures and the other one on classical control systems theory. Sometimes it is not possible to make a rigid classification because there are methods that combine techniques of one or the other alternative. In fact, problems of system behavior variation over time were considered by contributions from the area of adaptive control. In this sense, minimum variance controllers were designed for control loops with infrequent measure of the output variable in the case of fi order systems with delay [12]. These controllers, performing with infrequent observations of the output, were examined making use of the deterministic method introduced by [13]. In both perspectives, there are problems that are added to the absence of certain data such as outliers (data located far from the rest of the data), delays, measurement noises, and the variation of the system's behavior. In this last sense, an inferential control loop can also be analyzed and designed with feedback control techniques.

From the computer control loop perspective, and for slowly sampled systems, since the end of the 1980s the multirate (MR) option has been considered. An MR system is one where several sampling frequencies remain. If there are only two sampling frequencies, the system is called dual-rate (DR). It is especially important the case with slow output and fast input called MRIC (multirate input control). In a DR system is usual to consider an integer relation between the sampling periods and a regular pattern of sampled signals without sampling time mismatch. Actually, in this case, the control system is able to manage a slow measurement sampling frequency -far from an ideal sampling frequency, which can be chosen in order to achieve some performance-, considering an N times faster frequency, which is close to the ideal one, as the actuator frequency. As it is known, it is possible to obtain a linear model from a DR system using the discrete lifting technique as it will be explained later. Although initially MR control was used in the cement industry [14], especially in kilns [15] due to the diff t frequencies of chemical processes and material transport, and also in some heat exchanger processes [16], it was in the area of mechanical motion control [17] and robotic manipulators [18] where a higher number of applications were developed. In the DR case there were some attempts to adopt the adaptive control tendency as in [19] and as in contributions of relevant authors of inferential control [20, 21]. In recent years it has become a technique widely used in precision robotic manipulators [22], hard disk drive control [23], networked control [24, 25], industrial process [26], in UxV applications [27, 28] and in mobile robots [29, 24].

The starting point of this work is the simple dual-rate (DR) inferential control strategy introduced in [30]. In that case, with the original control loop including the original digital controller, the lost samples are provided by a process model. Only one of each $N$ real closed loop output samples can be measured. It is a DR control problem that can be approached by other strategies such as designing a nonconventional controller (slow input-fast output) computed from the process model and the original digital controller [31]. Both are time-varying discrete-time systems (TVDT). The current paper compares these two strategies using frequency response (FR) analysis techniques in order to reach a conclusion about the tracking and disturbance robustness. Control techniques based on FR can be adapted to deal with the uncertainty of the process control. QFT is an engineering technique for analysis and design of uncertain feedback systems that uses frequency domain specifications. QFT can be considered as a framework to design robust controllers for processes in which the uncertainty is typically parametric, that is, the process represented by a continuous transfer function $P(s)$ belongs to a family of plants $\mathbf{P}$ whose parameters values are included in certain fi intervals. The FR of all this family is called a "template". The specifications are given in the frequency domain in terms of admissible "bounds" on the FR of the closed loop transfer functions among the diff t selected inputs and outputs of the closed loop in order to achieve robust performance and stability specifications using a controller $G_R(s)$. For a number of frequencies, these specifications are combined with the system uncertainty description for obtaining constraints usually referred to as "boundaries" in the QFT procedure. For each of the selected frequencies, boundaries are given as curves in the Nichols plane, delimiting allowable regions for the open-loop transfer function $G_R(jw)P_0(jw)$ being $P_0$ the nominal process which can be any element in $\mathbf{P}$. Sometimes, an additional precompensator is needed. For a complete introduction to QFT see [32, 33]. QFT has been successfully applied to a wide variety of practical engineering problems: stable and unstable systems [34], single-input single-output and multiple-input multiple-output systems [35], linear and non-linear processes [36] and so on.

Therefore, the contributions of this paper are:

- It proposes two simple and practical DR control strategies to deal with vacant or missing samples (lost by any failure or deliberately lost) in a single-rate control.

- It provides an original QFT-based analysis procedure, never applied to DR systems to the best of the authors' knowledge, to study the effect of mismatch with respect to the real plant (stability and robustness). A new DR frequency response method based on discrete lifting modeling from [37] is used for this purpose. This method provides a tool to decide the best option in a specific case.

- It is concluded that the MBDR option outperforms the IC one when there are measurement losses, disturbances and the model of the process contains uncertainties.

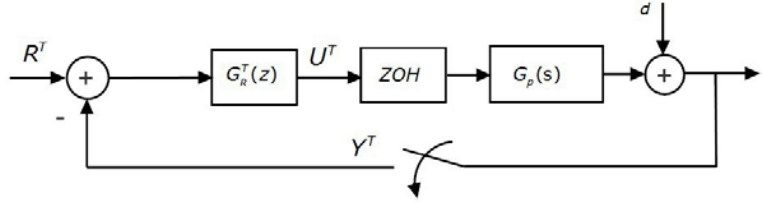

Figure 1: Fast digital control

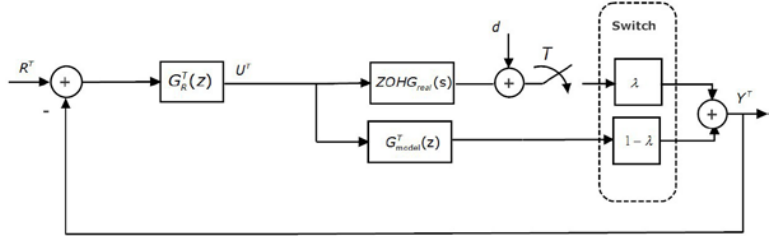

Figure 2: Dual-rate inferential control

## 2. Problem Statement

The problem to solve is based on a classical digital control with sampling period $T$. From a continuous-time process model $G_p(s)$ and a desired performance, a digital controller $G_R^T(z)$ is designed by any suitable direct or indirect method (where variable $z$ is used for the linear time-invariant (LTI)-transform argument at sampling period $T$). Under these conditions the control strategy is shown in Figure 1. In order to compare it with the other options in this paper, this control will be called "fast" control. $Y^T$ is used to describe the z-transform of the sequence $\{y(kT)\}$ [1].

When it is not possible to have access to all output samples, then diff rent options can be considered. Basically two options will be analyzed in this contribution: the fi one is what is known as Dual-rate inferential control (IC) and is shown in Figure 2. The second one is known as classical Model-based dual-rate control (MBDR) and will be explained later.

In Figure 2, the sampling of the real output provided by a $T$ discrete-time system or a continuous system sampled at period $T$ is known each $NT$ time instants, being $N \in \mathbb{N}$. This situation has been described using a variable $\lambda$

---

[1] Similarly $Y^{NT}$ for $\{y(kNT)\}$ will be used. This notation will be used for either signals or systems.

defined as:

$$\lambda = \begin{cases} 1 & \text{if } t = kNT \\ 0 & \text{if } t \ne kNT \end{cases} \quad k \in \mathbb{Z} \tag{1}$$

$\lambda$ represents a kind of switch that allows to get real samples if $\lambda = 1$ or samples each time instants $T$ delivered by a model $G^T_{model}(z) = Z_T[ZOHG_{model}(s)]$ in the other case [2]. Obviously a linear periodically time-varying system is met. It is called DR in the sense that two frequencies appear at $1/T$ and $1/NT$.

A classical DR solution is introduced in Figure 3 where the output is sampled at $NT$ (described by $Y^{NT}$) but the control is updated every $T$ ($U^T$) that is, $N$ times faster. To handle this solution, a special controller can be designed. In [31], a special strategy is proposed. This nonconventional structure controller consists of a slow part $G_1^{NT}(z^N)$, a digital hold that repeats the expanded slow controller output $N$ times,

$$H^{NT,T}(z) = \frac{1 - z^{-N}}{1 - z^{-1}}$$

and a fast controller $G_2^T(z)$. The slow and fast parts of this DR controller can be designed in different ways. The simplest form is to consider $G_1^{NT} = 1$ and $G_2^T = G_R^T(z)$, with unpredictable behavior. [38] describes formal design procedures. In this contribution, the MBDR controller will be used [31]. Basically, if $M(s)$ represents the desired closed loop performance of the original continuous system design

$$M(s) = G_R(s)G_p(s)/(1 + G_R(s)G_p(s)),$$

the MBDR will be designed with

$$G_1^{NT}(z^N) = \frac{1}{1 - M^{NT}(z^N)},$$

and

$$G_2^T = M^T(z)/G_p^T(z)$$

where $M^T(z)$ and $M^{NT}(z^N)$ are the $M(s)$ discretized with ZOH for periods $T$ and $NT$ respectively, and $G_p^T(z)$ is the continuous process model discretized with ZOH for period $T$. This procedure leads to a behavior that perfectly matches $M(s)$ at the sample points $\{y(kNT)\}$ but sometimes may introduce a ripple between samples due to the magnitude of $NT$. The way to overcome this ripple, if it appears, is as described in [31], considering

$$G_2^T = G_R^T(z)/(1 + G_R^T(z)G_p^T(z))$$

being $G_R^T(z)$ a proper discretization of $G_R(s)$ at period $T$. Therefore, the MBDR controller does not cancel the numerator of the process transfer function avoiding the ripple.

---

[2]ZOH is used to denote a zero order hold

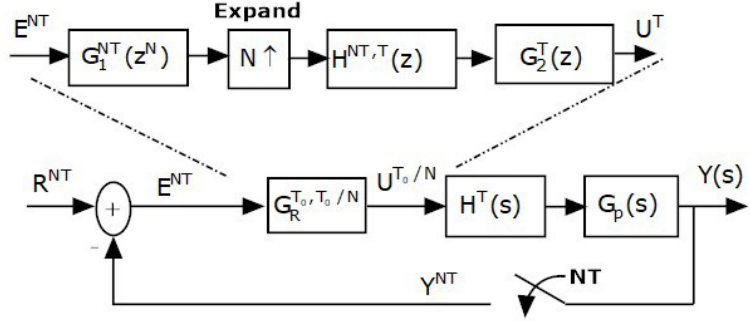

Figure 3: Model-based dual-rate control

3. **Discrete Lifting Modeling**

Some alternatives will now be considered: a fast $T$ digital control, an inferential DR control, and a classical DR control strategy. The last two options introduce a linear periodically TVDT system. Therefore, a procedure to transform it into an LTI system is convenient to make its treatment easier. A well-known method is the discrete lifting modeling. The original idea was introduced by Kranc [39] and it was called vector switch decomposition (VSD). In Figure 4, an open loop DR system is shown in which the input and output sequences have diff t sampling periods, $T_u$ and $T_y$. It is considered that they are rationally related. In Figure 5, the application of VSD to this case is shown considering the existence of integers $N_u$, $N_y$ such that $T_0 = T_u N_u = T_y N_y$ (indeed, then $T_u/T_y = N_y/N_u$ is a rational number) being $T_0 = lcm(T_u, T_y)$, where $lcm$ is the least common multiple. This metaperiod $T_0$ is the repetition period of every sequence in this scenario. With these conditions it is easy to introduce the discrete lifting operator that will allow an MR control strategy to be transformed into a $T_0$ single rate problem. This technique was explained with diff t details in excellent contributions like [40, 41, 42].

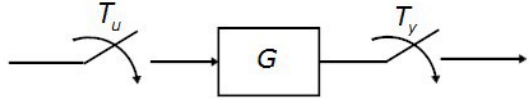

Figure 4: Generalized dual-rate system

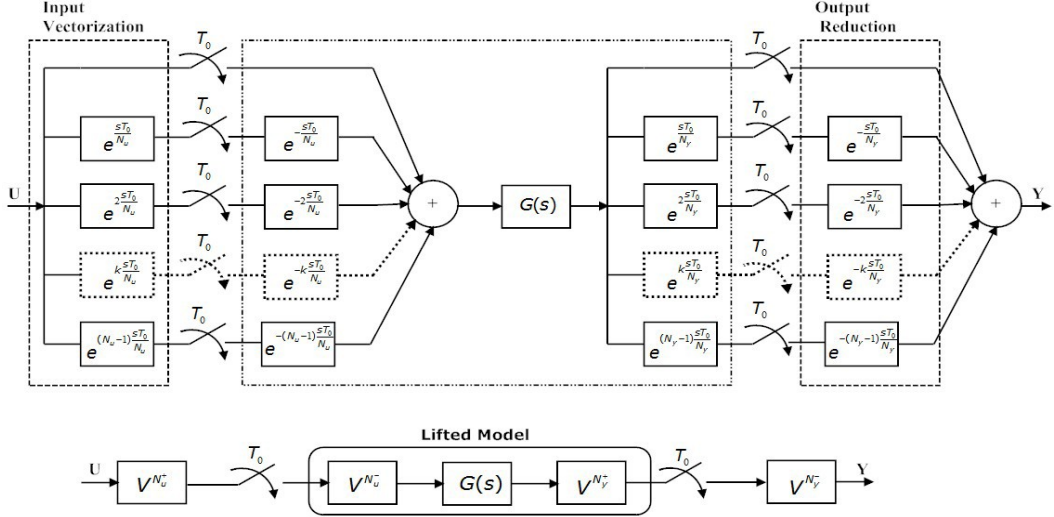

Figure 5: Kranc vector switch decomposition. General case

Kranc's basic idea is explained briefly with the intention to help understand discrete lifting through this perspective. A general sampler $T$ in a certain strategy such that $T_0 = NT$ ($T_0$ is the *lcm* of all samplers in the strategy) leads to a VSD with $N$ branches. An explanation of this case shows the concepts of vectorization and reduction. If $r(t)$ is a continuous signal, with Laplace transform $R(s)$, its sampling at period $T$ leads to the discrete sequence $\{r(kT)\}$ where the Z function will be noted by $R^T$. In these conditions, Kranc's idea was to propose a decomposition of this sampler into $N$ branches, each of them at period $T_0 = NT$

$$R^T = R^{T_0} + (e^{sT}R)^{T_0}e^{-sT} + \ldots + (e^{s(N-1)T}R)^{T_0}e^{-s(N-1)T}$$

So, it is possible to express

$$R^T = \begin{bmatrix} 1 & e^{-sT} & \ldots & e^{-s(N-1)T} \end{bmatrix} \left( \begin{bmatrix} 1 \\ e^{sT} \\ \vdots \\ e^{s(N-1)T} \end{bmatrix} R \right)^{T_0} \quad (2)$$

As seen in Figure 5, being in this case $N_u = N_y = N$ and $T_u = T_y = T$, the following notation is introduced

$$V^{N+} = \begin{bmatrix} 1 \\ e^{sT} \\ \vdots \\ e^{s(N-1)T} \end{bmatrix}$$

for the vectorization operation [3], and

$$V^{N-} = \begin{bmatrix} 1 & e^{-sT} & \ldots & e^{-s(N-1)T} \end{bmatrix}$$

for the reduction operation.
This leads to the formula

$$R^T = V^{N-} \left( V^{N+} R \right)^{T_0}$$

which is the representation of the $T$ sampler decomposition into $N$ branches including a new $T_0 = NT$ (*lcm*) sampler. Therefore, with this procedure it is possible to adapt every sampler in a general strategy with regard to a least common multiple ($T_0$) period of all sampling periods, leading to a multivariable LTI $T_0$ single rate system.

If the general case is being treated, then

$$(1 \; z^{-N_y} \; \ldots z^{-(N_u-1)N_y}) \left\{ \left( \begin{bmatrix} 1 \\ e^{sN_y} \\ \vdots \\ e^{s(N_u-1)N_y} \end{bmatrix} R \right)^{T_0} \right\} = V_d^{N_u^-} (V^{N_u^+} R)^{T_0} \quad (3)$$

---

[3]The $e^{mTS}$ blocks are used for off-line analysis purposes. For more explanation see 43].

where $z$ stands for greatest common divisor ($gcd$) period $T$ and $T_0$ is the least common multiple ($lcm$) of all sampling periods involved in the general strategy. In equation (3), the subindex $d$ (discrete) represents the substitution $z = e^{Ts}$. According to this procedure, the strategy in Figure 5 can be described as

$$Y^{T_y} = V_d^{N_y^-} \tilde{G} \left( V^{N_u^+} U \right)^{T_0} \quad T_0$$

where $\tilde{G}$ is the $T_0$ lifted matrix. Note that every signal is lifted in a $T_0$ sampling frame (obviously considering its own sampling period). So, it is common to model the behavior of the DR system characterised via a "lifted" transfer function matrix

$$y_l(z_{T_0}) = \tilde{G}(z_{T_0}) u_l(z_{T_0}) \qquad (4)$$

where the subindex $l$ denotes "lifted" and $z_{T_0}$ is referred to the $z$ variable at $lcm$ period $T_0$. In equation (4), $y_l$ is a vector of length $N_y$, $u_l$ is a vector of length $N_u$ and $\tilde{G}$ is a $N_y \times N_u$ transfer function matrix [40]. The lengths of the vectors are increased in the case of multiple-input and multiple-output (MIMO) systems (multiplied by the number of outputs and inputs, respectively). In Figure 5, if the block $G$ between samplers (it could be a continuous process plus zero-order hold) is a single-input single-output LTI system, the system in Figure 4 is clearly a shift-varying discrete-time system. In [41] it is proved that the $T_0$ lifted matrix $\tilde{G} = V_d^{N_y^+} G V_d^{N_u^-}$ is LTI $N_u$ input $N_y$ output shift-invariant. This correspondence between periodic and expanded LTI systems preserves both the analytic and algebraic properties of the systems [44]. In particular, block diagram algebra procedures could be considered [4]. In the next section, this algebra will be introduced to model the periodic systems to be compared.

Note that in Figure 5, $T_0/N_u = N_y T$ and $T_0/N_y = N_u T$. Using Mason's rule applied to this strategy, it is possible to prove that the central part of the

---

[4] It must be noted that it is also possible to define the lifted strategy from a discrete system, but obviously the input sampling period will be the same as the one applied to the system

diagram would be

$$\begin{pmatrix} 1 & z^{-N_u} & \ldots & z^{-(N_y-1)N_u} \end{pmatrix} \frac{\tilde{G}(z^{N_y N_u})}{N_y} \begin{pmatrix} 1 \\ z^{N_y} \\ \vdots \\ z^{(N_u-1)N_y} \end{pmatrix} \quad (5)$$

At this point, it is important to point out that the discrete-lifted operators can be expressed in an external [41] or internal [45] representation. In the contribution [42], the links between both were described. In the next subsection, two different examples are developed in detail.

### 3.1. Examples of discrete lifted systems

The procedure to obtain a lifted model of a system is explained in this subsection. Two cases are considered in this contribution: a continuous-time system preceded by a zero-order hold and a pure discrete system performing at the input sampling period but with a different output sampling period. In the first case, it is usual to consider a discrete state space realization $(A, B, C, 0)$ [5] at gcd period $T$ as a starting point. The objective is to obtain a realization $(A_l, B_l, C_l, D_l)$ at lcm $T_0$. These matrices are obtained, as it was mentioned before, by successive substitutions of the equations at sampling period $T$. As an example

$$\begin{aligned} y(kT_0+\eta T) &= Cx(kT_0+\eta T) = \\ &= C[A^\eta x(kT_0) + A^{\eta-1}Bu(kT_0) + \\ &\quad + A^{\eta-2}Bu(kT_0+T) + \ldots + Bu(kT_0+(\eta-1)T)] \end{aligned} \quad (6)$$

for $\eta = 1, \ldots, (N_u-1)N_y$. The ZOH implies that

$$u(kT_0 + pN_y T) = u(kT_0 + (pN_y+1)T) = \cdots = u[kT_0 + ((p+1)N_y - 1)T]$$
$$\forall p = 0, 1 \ldots, (N_u - 1) \quad (7)$$

In the second case, a discrete system with different output sampling period,

---
[5] A strictly proper process is considered.

equation (6) is valid but now

$$u(kT_0 + pN_yT) \text{ has a value}$$
$$u(kT_0 + (pN_y + 1)T) = \cdots = u[kT_0 + ((p+1)N_y - 1)T] = 0 \quad (8)$$
$$\forall p = 0, 1, \ldots, (N_u - 1)$$

is fulfilled.

The matrices of the lifted realization are obtained, in each case, by appropriately stacking the terms from the above equation; the lifted matrix $\tilde{G}$ is derived from this quadruple [42].

## 4. **Dual-rate systems frequency domain**

A description of the frequency domain procedures used for the analysis comparison of both strategies is introduced in this section. First, a new efficient algorithm [37] for computing the DR systems frequency response is referenced and how to understand this algorithm is also described. The objective is to analyze the robust stability and robust tracking of diff t systems with the options considered in this work. The option selected for a unifi analysis is the quantitative feedback theory (QFT) [32, 46].

### 4.1. DR Systems Frequency Response

Theorem 1 in [37] establishes that it is possible to calculate the frequency response of a DR system once the LTI representation of the DR system (the lifted model) has been obtained. In a practical way, considering (5), the frequency response computation is given by the product of the frequency response of a left factor $[1 \ z^{-1} \ z^{-2} \ \ldots \ z^{-(N_y-1)}] G_l(z^{N_y})$ replacing $z = e^{j\omega_r T_y}$, which gives a row vector, and the right factor (column vector) $(1 \ z \ z^2 \ \ldots z^{N_u-1})^T$ replacing $z = e^{j\omega T_u}$.

[47] explains that for an input frequency $b$ Rad/s, the output is the sum of sinusoidal signals with frequencies $b$, $b+N_u\omega_s$, $b+2N_u\omega_s$, ... $b+(\frac{N_y}{gcd(N_u,N_y)} - 1)N_u\omega_s$ with amplitude and phase determined by the Bode diagrams at the indicated points (note that $\omega_s = 2\pi/T_0$). Therefore, in the coprime case (assumed in this contribution), $N_y$ components are obtained. Using the $T_0$ resampling of each of these $T_y$ components, a $T_0$ pure sinusoidal is achieved (for details see [47]).

## 4.2. Dual-rate systems quantitative feedback theory

Obviously, this work does not try to show or revisit the quantitative feedback theory (QFT) [32, 46]. The purpose is just to introduce its application to a DR system. The starting point is [48] that introduces the QFT for direct digital control.

In the present work, the QFT analysis methodology is slightly diff rent than usual. Classically, QFT analysis and design involves the mapping of robust specifications (robust stability, robust disturbance rejection, control effort, robust tracking among others) into certain bounds on a nominal loop transmission followed by loop shaping for some selected frequencies. These bounds split the complex plane into two areas where the open loop transfer function should lie inside one. In this work QFT is used as an analysis tool. The controller design (IC and MBDR) is fi      performed considering a nominal plant model and design specifications (relative stability with some stability margins, disturbance rejection, ...). For this nominal design, QFT is used for computing restrictions (boundaries) to be satisfied by the nominal open-loop gain, for diff ent working frequencies (note that templates are reduced to a point in the Nichols plane). Once the nominal controller design is performed, the effect of plant uncertainty on the satisfaction of (closed loop) design specifications is considered. This is done by computing the open-loop gain corresponding to the plant with uncertainty, and checking if its FR satisfi  the previously computed boundaries. For all these calculations, the DR frequency response introduced before and the algebra with lifted blocks that will be explained in section 5 are absolutely necessary.

In this work two of the specifications are considered but obviously the procedure can be applied for all the options. Now the gain (GM) and phase margins (PM), and the discrete-time output disturbance rejection are considered. In both cases and with respect to a structure like the one shown in Figure 1, some restriction must be planned in the closed loop transfer function denominator [48]. In the case of robust stability margins $|1 + G_p(z)G_R(z)| \geq \mu$ with $z = e^{jwT}$, $w \in [0, w_s/2]$ being $GM = (1 - \mu)^{-1}$ and $PM = 180° + 2cos^{-1}(\mu/2)$. For the output disturbance rejection can be considered $|1 + G_p(w)G_R(w)| \geq (1/\delta(w))$. In the applications, sections 6 and 7, these restrictions will be considered more specifically.

## 5. Closed loop lifted models

In this section, the block diagram (obtained by lifting) algebra is used to show the models of the closed loops to be compared. Obviously, it is not necessary to obtain the fast digital control lifted model. The closed loop $Y$ vs $R$ for MBDR and IC cases are shown in Figures 6 and 7 [6]. These fi    show the complete evolution from the initial digital control strategy to the equivalent one using lifted signals. The substitution of each sampler following Kranc's idea (vector switch decomposition) is the fi step. A reference sampled at whatever $T_{ref}$ and an output sampled at $T_3 = T_0/N_3$ with $N_3$ high enough have been added to the original strategies. Note that the sampling for feedback purposes is only possible at $T_1 = T_0/N_1$, that is, a low magnitude. $T_3$ reflects open-loop intersampling. It is not realistic, but it is included with the intention to define a general problem. As it can be seen due to the reasons introduced in section 3.1, in a middle stage *kranczoh* or *krancdig* are diff  tiated corresponding to lifting matrices of ZOH plus continuous process and digital controller respectively. Then the notation $K_{type}G^{N_u,N_y}$ is considered where the superindex is referred to the multiplicity of input and output regarding the *lcm* $T_0$ period, and the subindex is related to the kind of operator. In this notation, the transfer function $G$ is also added for easy understanding. In Figures 6 and 7, the input control is updated every $T_2 = T_0/N_2$. Note that $N_3 \geq N_2 \geq N_1$. Analyzing these fi     and the procedure explained graphically, it is possible to express the lifted input-ouput relation for the MBDR loop (9) and for the IC Loop (10). Note that in this last case, IC, the model has been lifted at the metaperiod. In this last case, the switch selector which models the lost information and the complementary information delivered by the model leads to a system with periodicity $T_0$.

$$\bar{U} = K_{DIG}G_2^{N_2,N_2} \times K_{DIG}G_1^{N_1,N_2} \times \left[I_{\alpha_1 \times \alpha_2} + K_{ZOH}G_{preal}^{N_2,N_1} \times K_{DIG}G_2^{N_2,N_2} \times K_{DIG}G_1^{N_1,N_2}\right]^{-1}$$
$$Y^{T/N_3} = K_{ZOH}G_{preal}^{N_2,N_3} \times \bar{U}$$
$$Y^{T/N_1} = K_{ZOH}G_{preal}^{N_2,N_1} \times \bar{U}$$
$$(9)$$

---

[6] Note that in these figures the pure zero-order hold is shown as $ZOHG_{unit}(s)$ being $G_{unit}(s) = 1$

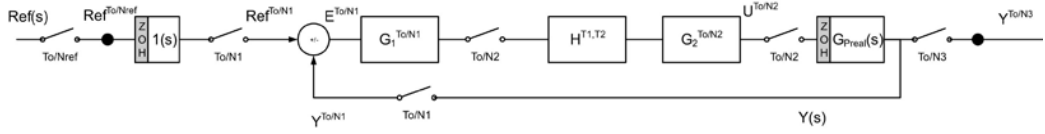

(a) Original diagram

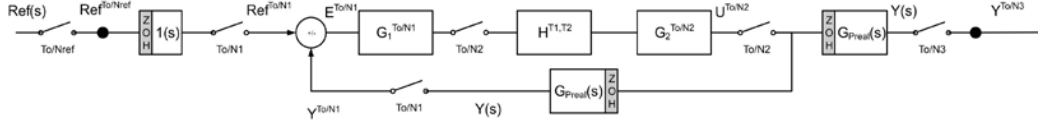

(b)

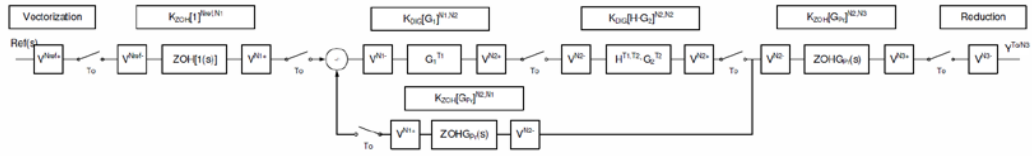

(c)

Figure 6: Y/R MBDR lifted process

where $\alpha_1 \times \alpha_2$ is the size of the product $K_{ZOH}G_p^{N_2,N_1} \times K_{DIG}G_2^{N_2,N_2} \times K_{DIG}G_1^{N_1,N_2}$ and the notation $I$ is used for the identity matrix.

For the case of the IC strategy

$$\begin{aligned}
\bar{U}_R &= K_{DIG}G_R^{N_2,N_2} \times \left[I_{\beta_1 \times \beta_2} + \tilde{K}_{ZOH}G_p^{N_2,N_2} \times K_{DIG}G_R^{N_2,N_2}\right]^{-1} \\
Y^{T/N_3} &= K_{ZOH}G_{preal}^{N_2,N_3} \times \bar{U}_R \\
Y^{T/N_1} &= K_{ZOH}G_{preal}^{N_2,N_1} \times \bar{U}_R
\end{aligned} \quad (10)$$

being $\beta_1 \times \beta_2$ the size of the product $\tilde{K}_{ZOH}G_p^{N_2,N_2} \times K_{DIG}G_R^{N_2,N_2}$

Now, $\tilde{K}_{ZOH}G_p^{N_2,N_2}$ is quite special because it is made up of two diff rent parts which correspond to the order of sampling in the switch of Figure 2. This will be easier to understand through an example.

### 5.1. Example

Take, for instance, $N = 3$ and consider that only the first of a set of $N$ samples is picked up from the real process $G_{preal}$ and the other $N - 1$ are delivered by the model $G_{pmodel}$. If both are second-order systems and

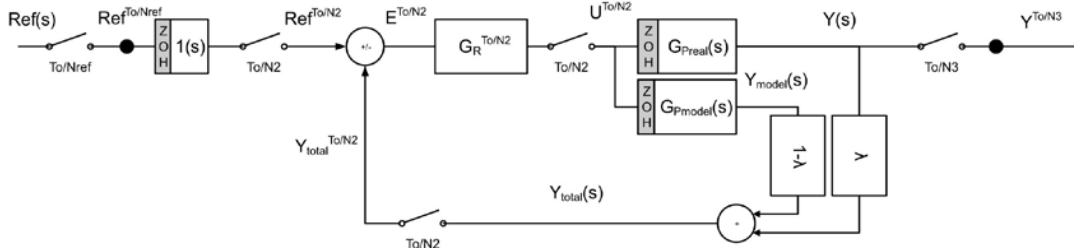
(a) Original diagram

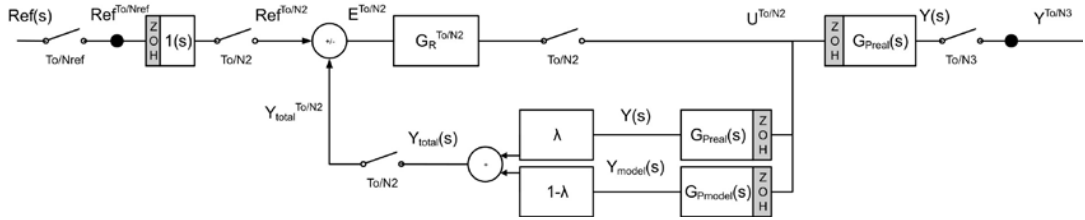
(b)

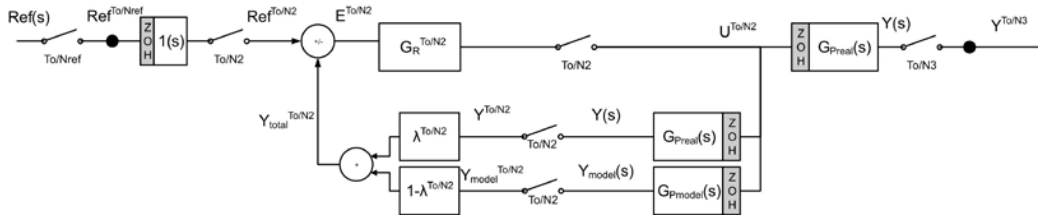
(c)

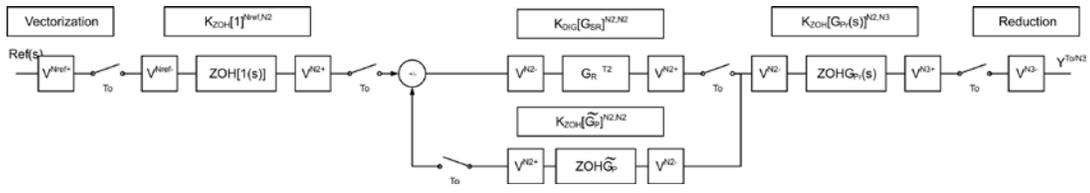
(d)

Figure 7: Y/R IC lifted process

described by a state representation $A_r, B_r, C_r, D_r$ and $A_m, B_m, C_m, D_m$, then $\tilde{K}_{ZOH}G_p^{N,N} \equiv A_t, B_t, C_t, D_t$. Some auxiliary matrices are defined to simplify the notation [7]

$$\chi_1 = \left[\begin{array}{c|c} I_{dimA_r} & 0 \\ \hline 0 & 0_{dimA_m} \end{array}\right] \quad \chi_2 = \left[\begin{array}{c|c} 0_{dimA_r} & 0 \\ \hline 0 & I_{dimA_m} \end{array}\right]$$

$$\chi_3 = \left[\begin{array}{c|c} 1_{1\times 1} & 0_{1\times N-1} \\ \hline 0_{N-1\times 1} & 0_{N-1\times N-1} \end{array}\right] \quad \chi_4 = \left[\begin{array}{c|c} 0_{1\times 1} & 0_{1\times N-1} \\ \hline 0_{N-1\times 1} & 1_{N-1\times N-1} \end{array}\right]$$

$$\begin{aligned} A_t &= \left[\begin{array}{c} \chi_1 \times A_r + \chi_2 \times A_p \end{array}\right] \\ B_t &= \left[\begin{array}{c} B_r \\ \hline B_m \end{array}\right] \\ C_t &= \left[\begin{array}{c|c} \chi_3 \times C_r & \chi_4 \times C_m \end{array}\right] \\ D_t &= \left[\begin{array}{c} \chi_3 \times D_r \\ \hline \chi_4 \times D_m \end{array}\right] \end{aligned} \quad (11)$$

In this section, the complete procedure is applied to two diff t simulations with diff ent conclusions. The fi step is to obtain the lifting modeling of MBDR and IC cases shown in section 3 for each example. Then, the algebra introduced in section 5 and the DR systems FR introduced in section 4, allows the consideration of the QFT procedure in order to analyze both cases.

---

[7]a $N, N$ is used in order to simplify the explanation. That is, $N_2 = N_3 = N$ which is not relevant for the analysis procedure.

*Example 1*

The first is a continuous system with model

$$G_{pmodel} = \frac{1.5}{(s+0.5)(s+1.5)}$$

this plant is controlled in a closed loop by a serial PID controller,

$$u(t) = K_p \left[ e(t) + T_D \frac{d}{dt} e(t) + \frac{1}{T_i} \int_0^t e(\tau) d\tau \right]$$

with $K_p = 8$, $T_D = 0.2$, and $T_i = 3.2$ in order to achieve some specifications. A discrete-time PID controller approximation for a sampling period $T$, is given by [49]

$$q_0 = K_p \left( 1 + \frac{T_D}{T} \right)$$

$$q_1 = -K_p \left( 1 + 2\frac{T_D}{T} - \frac{T}{T_i} \right)$$

$$q_2 = K_p \frac{T_D}{T}$$

Nevertheless, multiplicative uncertainty has been considered in the real system

$$G_{preal} = \frac{1.5}{(s+0.5)(s+1.5)} \frac{10s+1}{100s+1}$$

that is, with a low frequency pole and high frequency zero. A single-rate control, which is called "fast" is designed for $T = 0.1$ s. Due to diverse restrictions it is only possible to sample the output measurement every $T_0 = 0.3$ s. With these conditions the MBDR and IC with $N = 3$ is planned. Now the three options are going to be compared. In Figure 10, the closed loop step response is shown for the three controllers designed for the model plant (obviously the fast and inferential cases have the same response if ideal conditions are considered). The application to the real process (that is, with model plant mismatching MPM) has been considered in Figure 11. As it can be seen, the IC leads to a very slow output response, being the MBDR control response similar to the one of the fast single rate control. The closed loop DR Bode diagram in Figure 12 confirms this behavior. As it can be seen, the IC case has a small bandwith. The QFT analysis described in subsection 4.2 with two reasonable restrictions for robust stability margins and output disturbance rejection $\mu = 0.3$ and $\delta(w) = w/1.5$ is shown in Figures 13 and 14 for MBDR and IC cases respectively. A small gain increase in MBDR leads to the correct performance, but for the IC case even a high gain increment is not able to ensure the disturbance rejection for small frequencies and makes the system unstable. In Figure 14 the dashed line corresponds to the open loop shape incremented by a gain of 300. Following [50] the system is clearly unstable. The analysis has been performed with the help of the QFT Toolbox [51] with scripts adapted to the current problems.

Moreover, if the closed output response is analysed with respect to a step in the reference and an output disturbance $d(t) = sin(0.04\,t)$, (see Figure 15), the ineffective disturbance rejection capability of the IC is clear, while the MBDR control is valid for this purpose.

*Example 2*

The second example is taken from [30]. In this case, the discrete model process[8] and the PI controller designed are

$$G_{pmodel}(z) = \frac{0.15}{z - 0.9}, \qquad G_R^T(z) = 3 + \frac{0.5}{1 - z^{-1}}$$

The actual plant is

$$G_{preal} = \frac{0.13}{z^3(z - 0.92)}$$

Now a similar reasoning as in the example before is run. The closed loop step responses for the three controllers designed for the model plant are shown in Figure 16. In this case, the three closed loops have a very similar output reponse. For this example the MBDR design has been adapted from [30]. Now the continuous closed loop transfer function $M(s)$ is not known. Only $M^T(z)$ at $T$ is known. A simple discrete to continuous procedure was used for obtaining the discrete $NT$ transfer function $M^{NT}(z^N)$. As in [30] $N = 4$ was considered. The application to the process with MPM is shown in Figure 17. In this case, the IC leads to the best closed loop output response (as it is said in [30]), being the single control unstable. The closed loop Bode diagrams are shown in Figure 18 and confirm this behavior. With the same tool used in example 1, the QFT analysis explains that in this case the IC is able to perform within the restrictions for stability $\mu = 0.3$ and output disturbance rejection $\delta(w) = w/1.5$ as depicted in Figure 19. In the MBDR case both restrictions are not satisfied. A slight increment of gain (about 1.3) would provide the output disturbance rejection compliance but the system would be unstable ($w = 1.7, 2$ R/s) as indicated by Figure 20.

---

[8]it is possible to consider $T = 1$. This is not a relevant fact for the current purposes

However, if the real plant were

$$G_{preal} = \frac{0.4}{z^3(z - 0.96)}$$

that is, with the same transport delay but with diff t dynamics than the model, the IC would lead to instability and the MBDR would perform an acceptable system output response.

In general when there are significant diff in gain and/or constant times, the MBDR alternative is much better than the IC option.

## 7. Experimental Application

At this point, the contents and procedures discussed before will be applied to a diff tial mobile robot path-tracking. The robot housing is formed by three plastic pieces constructed with a 3D printer. An Arduino Due with a 32-bit Atmel SAM3X8E ARM Cortex-M3 microcontroller was used. The robot has two wheels that are operated by two Pololu MG37D motors (12 volts) with encoders and two caster balls. A lithium battery with three cells and nominal voltage of 11.1 volts was utilized. An HC-SR04 ultrasonic sensor was used for measuring distances. In Figure 21 two pictures of the robot are shown.

This unmanned ground vehicle (UGV) should follow a predefi path. A pure pursuit tracking algorithm was considered. The position of the mobile robot is obtained using an energy consuming sensor and therefore restricts the time of autonomous use. It is proposed that this sensor could work less frequently but it is unfeasible in a classical single rate control due to the servomotors dynamics. In other words, a large sensing period but a short control update is needed in this scenario. Therefore, the problem introduced in this work reflects this situation. Both options, MBDR and IC, were used in this application. In both cases a $N$-slow measurement and fast control can be considered as previously described.

### 7.1. Path tracking

The mobile robot must track a prescribed path. The fi problem is to know in advance what is the next point to reach, that is the target point. There are diff t methods: follow the carrot, vector pursuit, pure pursuit and so on. In this work, a pure pursuit algorithm was used. Basically, the method tries to compute one arc between the current position and the target point of the predefi trajectory. For this purpose it is necessary to choose the next point to reach. It is usual to fi what is known as a look ahead distance so that the algorithm fi s it. In this application, the pure pursuit method was selected because the path followed by the robot is smoother. The pure pursuit method and concepts like look ahead distance can be consulted in known references [52, 53]. Once the future point has been computed, the diff tial mobile robot kinematics lead to the proper angular velocity of both wheels in order to get the target point, [54]. In our case, the linear speed was prefi in a constant value that makes reaching the next point possible. Therefore, at every period the path tracking algorithm will provide the angular speed reference for each wheel once the current position is detected. A closed loop for each of the wheels was considered considering a slow measurement of actuators encoders. It would have been possible to use the encoders measurements to determine the robot position by odometry but the precision is much better with the selected procedure.

### 7.2. Dual-rate controller

The IC and MBDR options are planned by considering a slow measurement, that is, considering less sensor measures. If the sensor is used less frequently, it will consume less energy and therefore, battery operation life will be greater. This slow frequency can be processed as samples losses with respect to a fast frequency as it was already explained. In the IC case, the samples that have not been sensed are provided by an approximate actuator model, while in the MBDR case only the measurements obtained are used by the control algorithm.

In the set-up used, both servomotors (actuators) are considered to have the same transfer function (actually there are very slight diff

$$G_p(s) = \frac{W(s)}{V(s)} = \frac{0.1276}{0.1235s + 1} \qquad (12)$$

where the output $W$ is in rad/s, and the input $V$ in volts. The ideal control in order to achieve some specifications leads to a classical PI controller

$$u(t) = 6\left[e(t) + \frac{1}{0.12}\int_0^t e(\tau)d\tau\right]$$

If it is discretized [49], the single-rate controller for $T = 100$ ms is

$$G_r^T(z) = \frac{6z - 1}{z - 1} \qquad (13)$$

It was considered that the measurement should be $T_0 = 300$ ms, that is $N = 3$ for proper energy saving. With these conditions, the slow and fast controllers in the MBDR loop were

$$G_1^{NT}(z^N) = \frac{z_s^2 - 0.23z_s + 0.014}{z_s^2 - 1.08z_s + 0.08} \qquad (14)$$

$$G_2^T(z) = \frac{6.576z^2 - 5.78z + 1.27}{z^2 - 0.976z + 0.24} \qquad (15)$$

where $z$ stands for $T = 0.1$, and $z^N = z_s$ for $NT = 3T = 0.3$. As it was explained, a digital hold is included between the slow and fast parts

$$H^{NT,T}(z) = \frac{1 - z^{-N}}{1 - z^{-1}}$$

In this case, in order to apply the considerations of this work, it is considered that the model is

$$G_m(s) = \frac{0.2}{0.2s + 1}$$

that is, there is a strong uncertainty in gain and time constant. The QFT analysis, performed in a completely similar way than in 6 with $\mu = 0.6$ and $\delta(w) \leq w/1.5$ for robust stability and output disturbance rejection, is depicted in Figure 22. Figure 23 shows the real trajectory followed by the UGV when the ideal, mismatch and mismatch with output disturbance are considered. As it can be seen, the MBDR option is much better than IC one when the UGV must track an "H"-shaped path. The output disturbance considered was $d(t) = 2sin(0.2\ t)$. In regards to energy saving, when tracking the "H" reference, increasing the value of $N$ from 1 to 3, the power consumption dropped by about 27% in our experiments. The power was calculated measuring intensity and voltage using a CompactRio (National Instruments) system. In order to analyze the infl nce of the non-linearities of the servomotors (dead zone and saturation) a simulation was carried out. The results showed in Figure 24 also confirm that the MBDR option is better than the IC one.

8. **Conclusions**

In this contribution diff t DR model-based control design methodologies have been compared. These methods are motivated by restrictions that prevent the consideration of an ideal sampling pattern. The use of DR techniques allows to model this kind of situations, and a new efficient tool for DR FR computation is crucial for analyzing and detecting advantages or disadvantages for these diff t alternatives. It is very important how high the uncertainty of the real process is, that is, the comparison is problem-dependent but it seems that in any case the inferential control is unable to reject disturbances. A new QFT analysis procedure for DR systems has been used. Even though every case needs a particular analysis, when the model contains important uncertainties, the MBDR is a better option. Avenues for further research like the extension of this method to non-linear plants and the consideration of less simplistic inferential control strategies must be considered in the future.

9. **Acknowledgements**

This work is funded in part by grant RTI2018-096590-B-100 from Spanish government, and by European Commission as part of Project H2020-SEC-2016-2017 - Topic: SEC-20-BES-2016 - Id: 740736 - "C2 Advanced Multi-domain Environment and Live Observation Technologies" (CAMELOT). Part WP5 supported by Tekever ASDS, Thales Research and Technology, Viasat Antenna Systems, Universitat Politècnica de València, Fundação da Faculdade de Ciências da Universidade de Lisboa, Ministério da Defesa Nacional - Marinha Portuguesa, Ministério da Administração Interna Guarda Nacional Republicana.

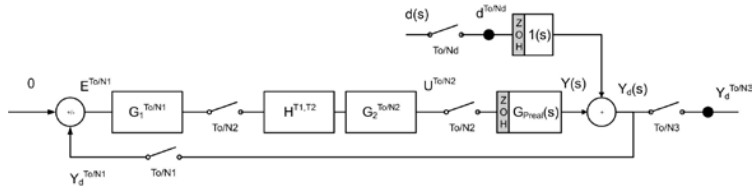

(a) Original diagram

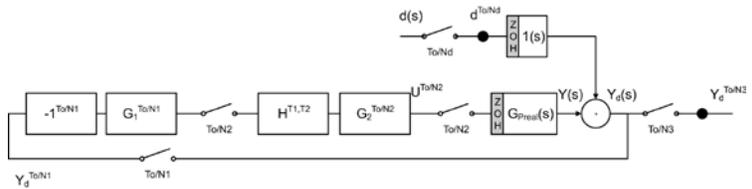

(b)

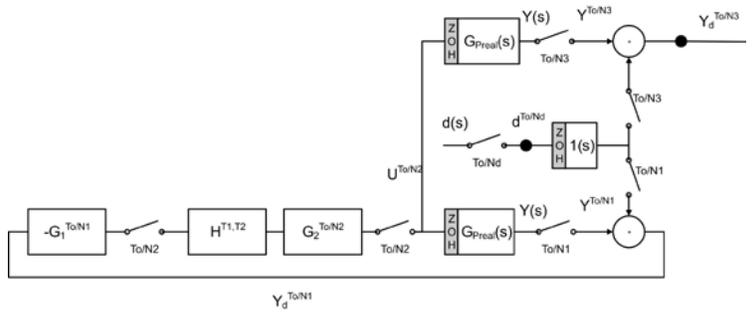

(c)

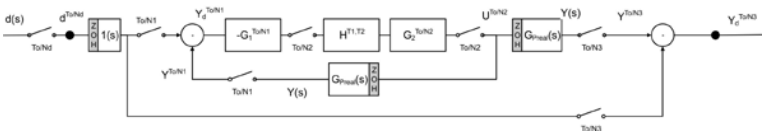

(d)

(e)

(f) Lifted diagram

Figure 8: Y/d MBDR lifted process

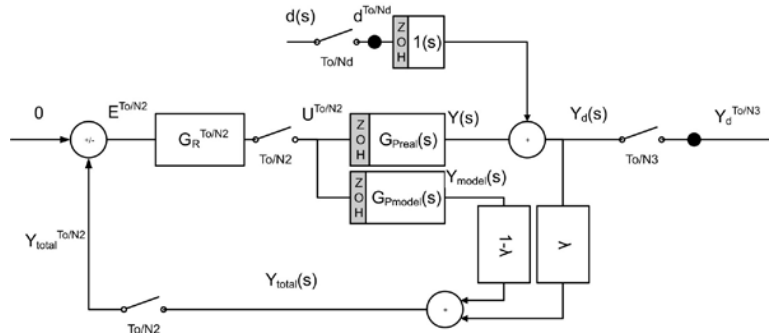

(a) Original diagram

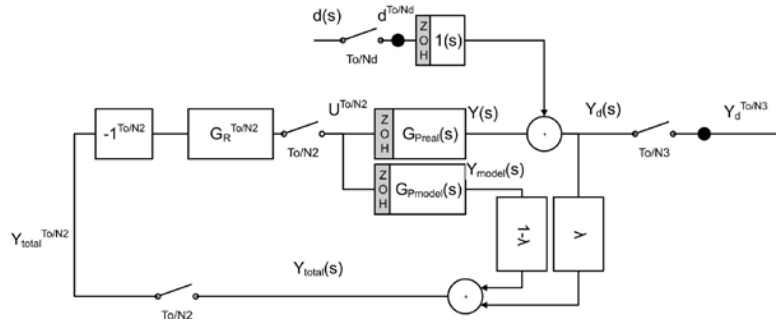

(b)

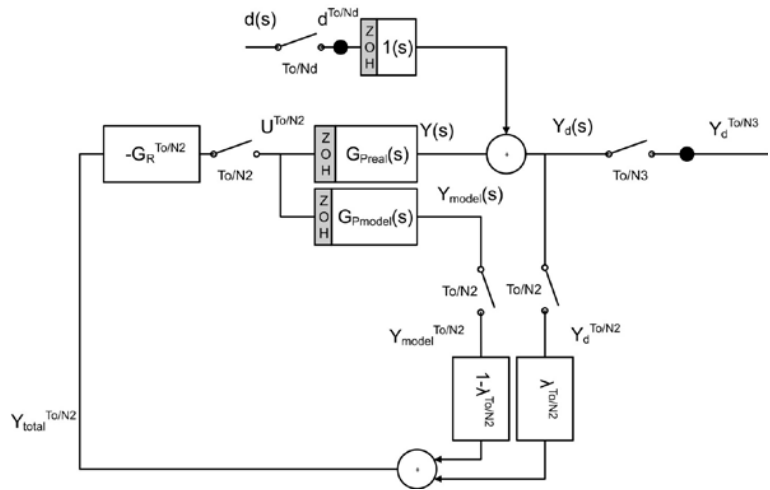

(c)

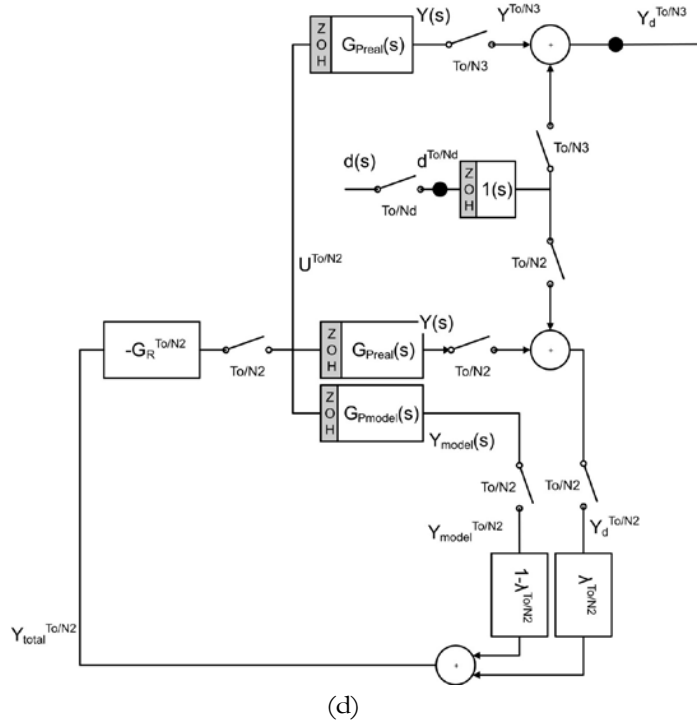

(d)

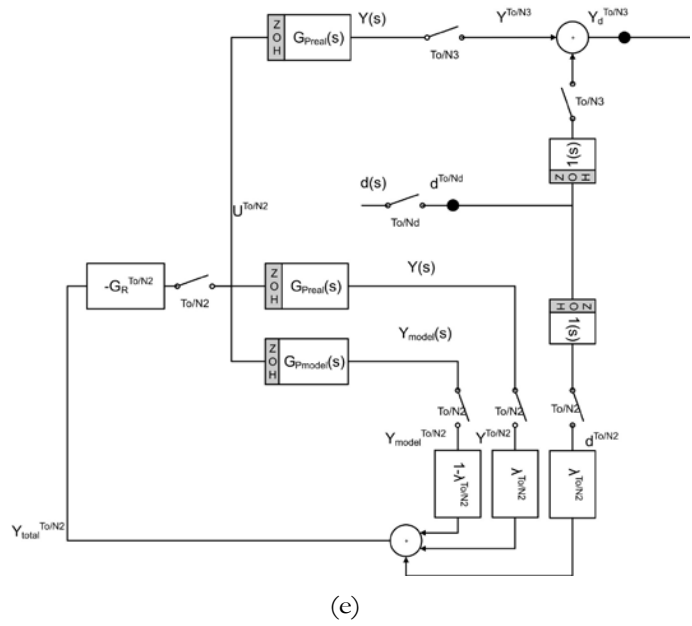

(e)

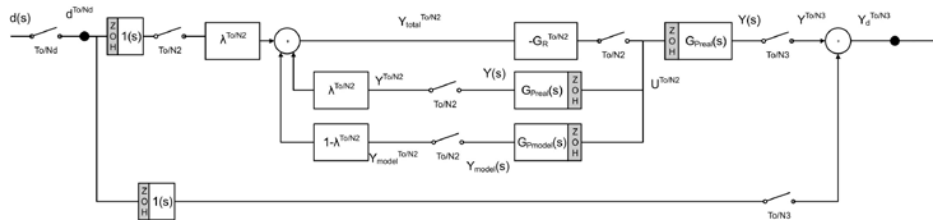

(f) Lifted diagram

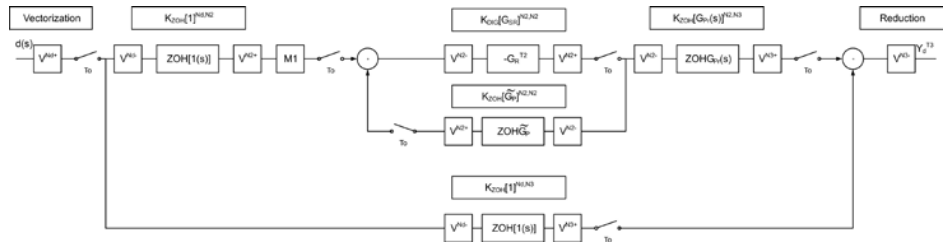

(g) Lifted diagram

Figure 9: Y/d IC lifted process

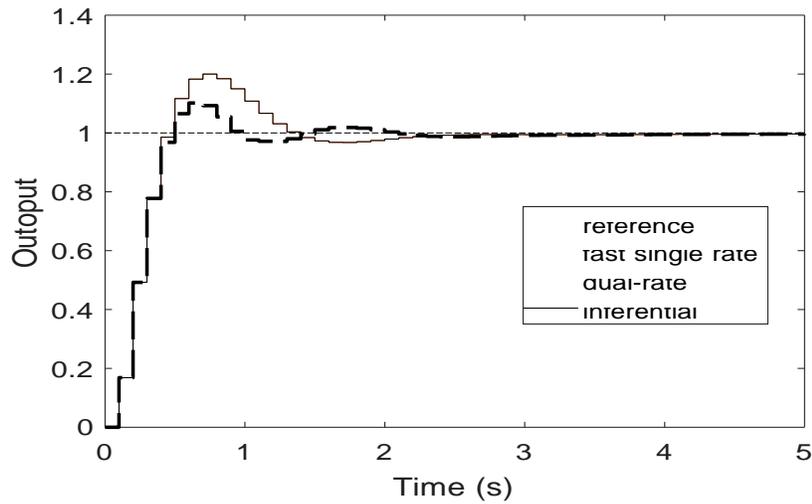

Figure 10: Example 1. Time response without MPM

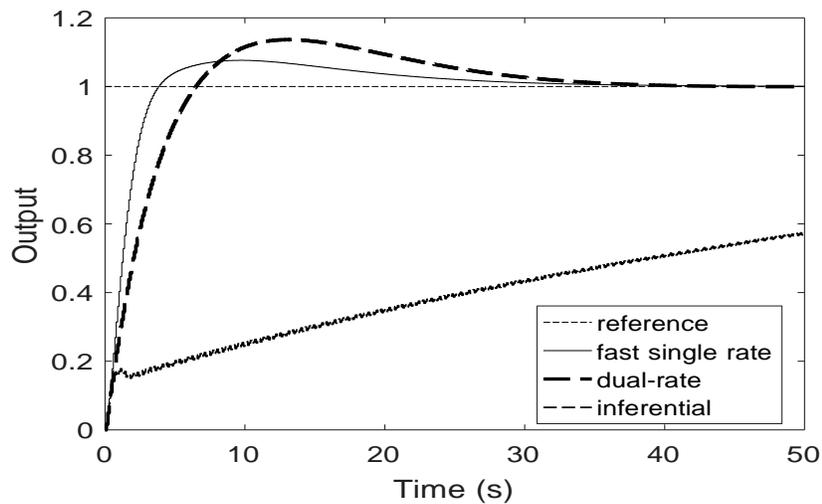

Figure 11: Example 1. Time response with MPM

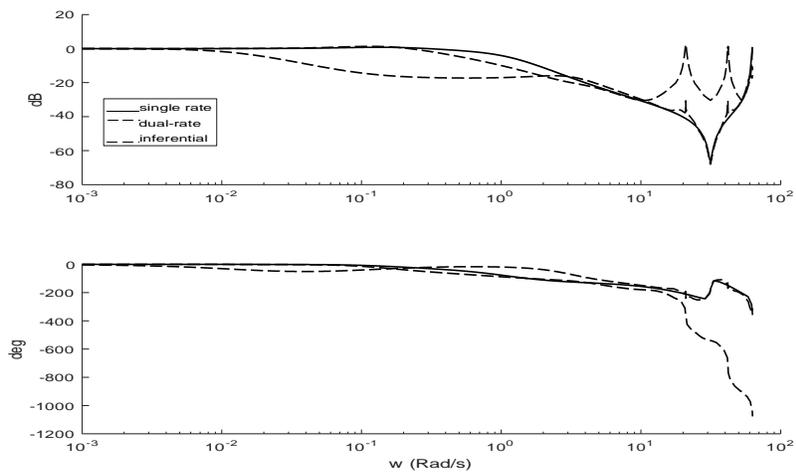

Figure 12: Example 1. Y/R DR systems Bode. Methods comparison

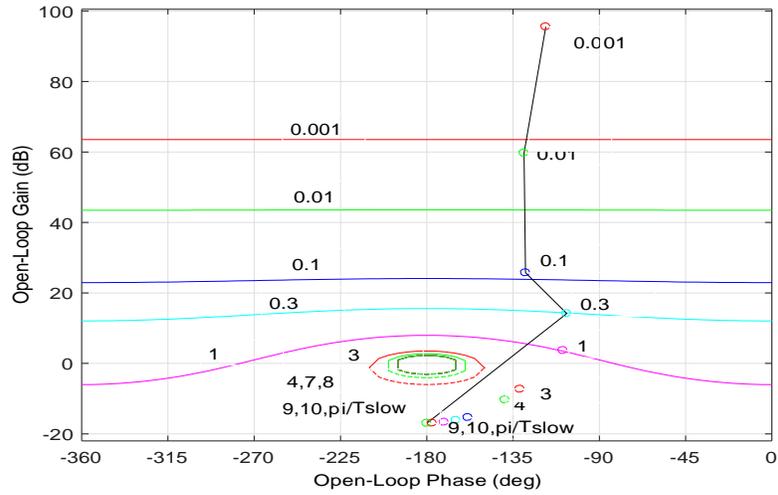

Figure 13: Example 1. MBDR case QFT analysis

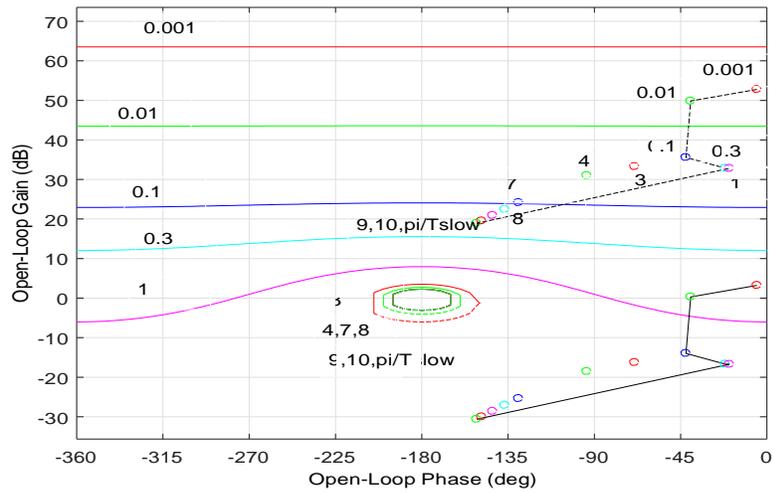

Figure 14: Example 1. IC case QFT analysis

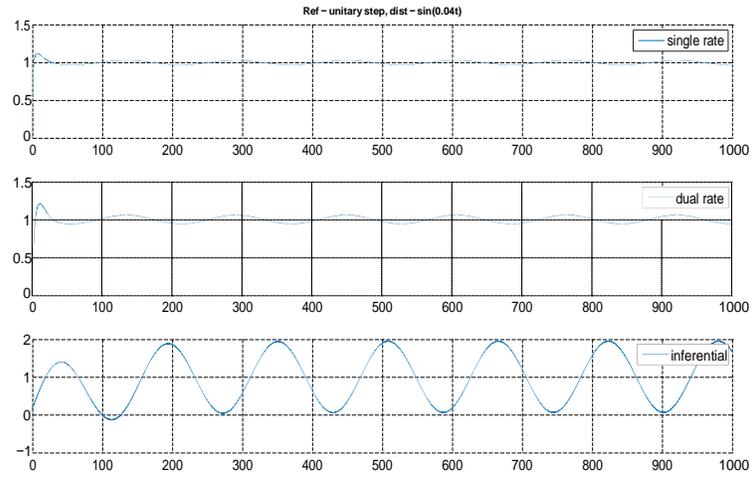

Figure 15: Disturbance rejection comparison

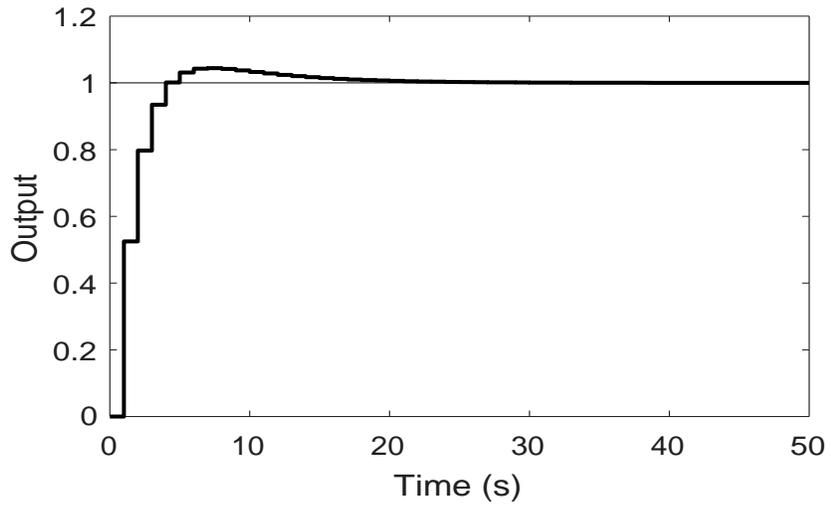

Figure 16: Example 2. Time response without MPM

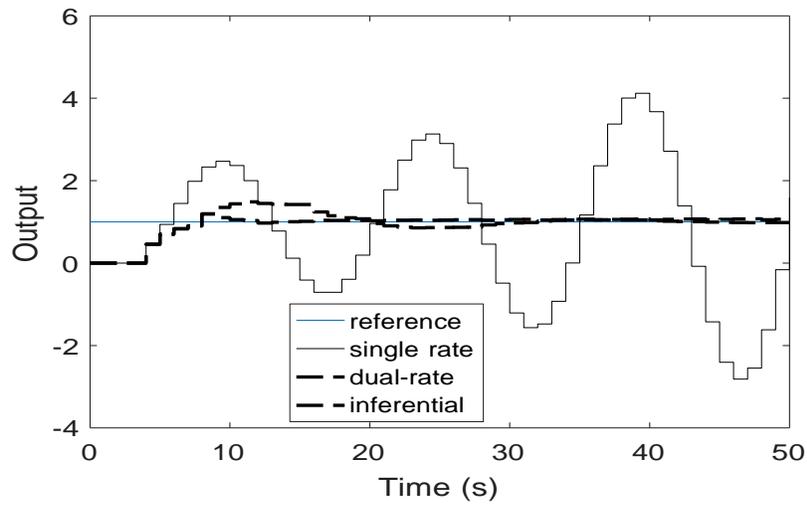

Figure 17: Example 2. Time response with MPM

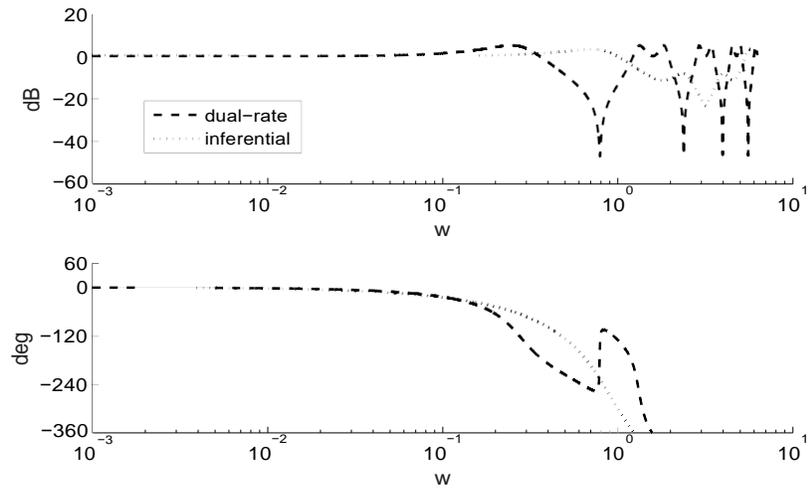

Figure 18: Example 2. Y/R DR systems Bode. Methods comparison

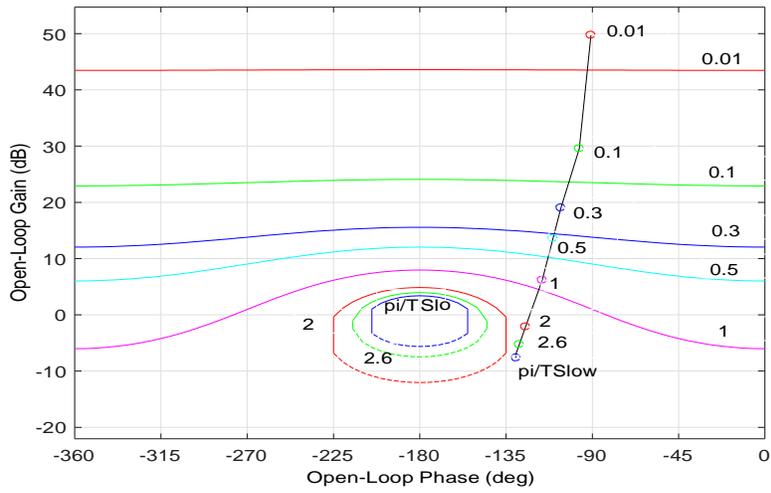

Figure 19: Example 2. IC case QFT analysis

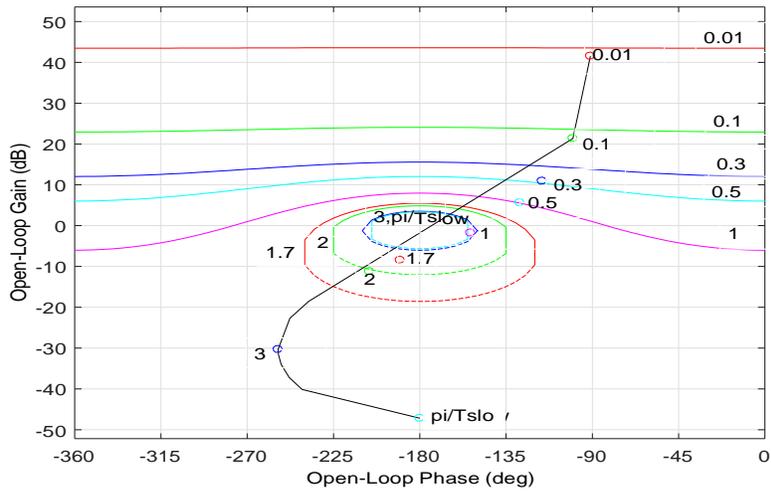

Figure 20: Example 2. MBDR case QFT analysis

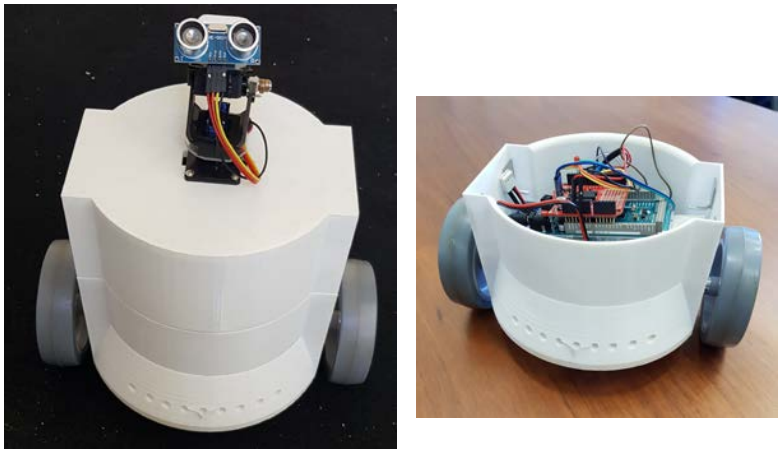

Figure 21: Autonomous robot used

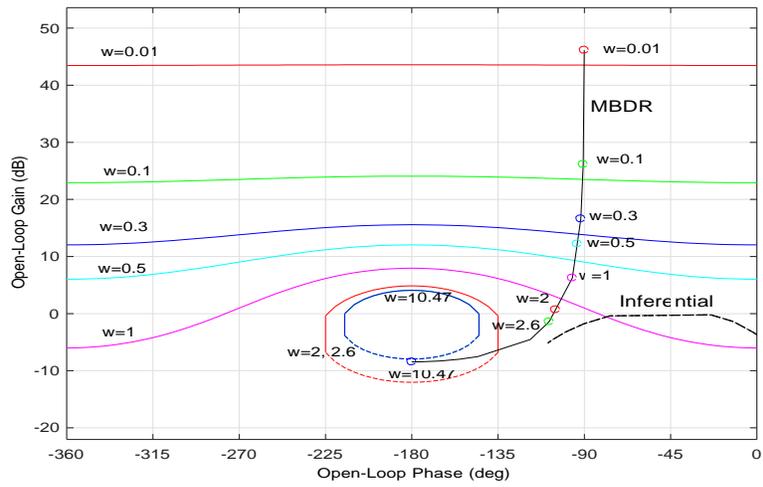

Figure 22: UGV QFT analysis

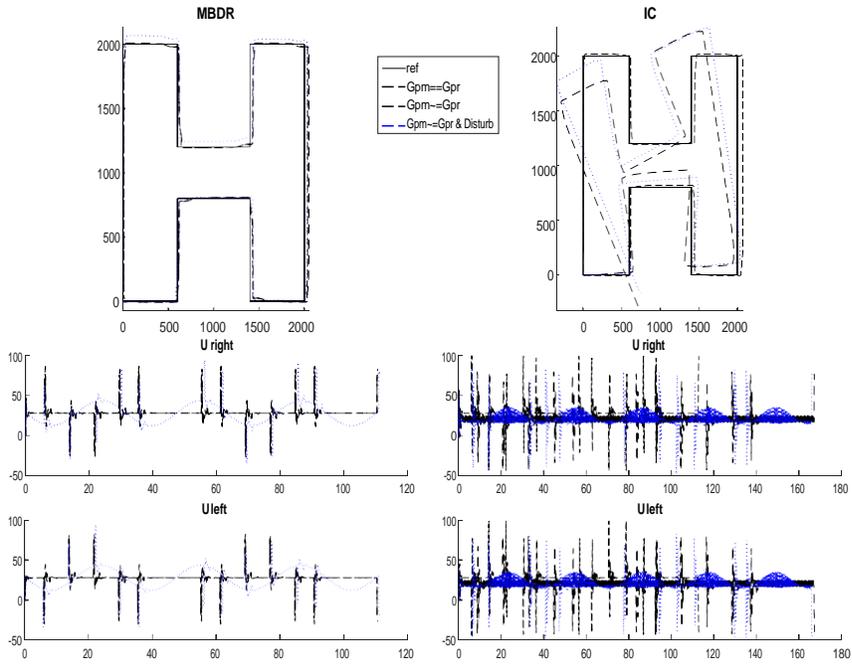

Figure 23: UGV Path tracking

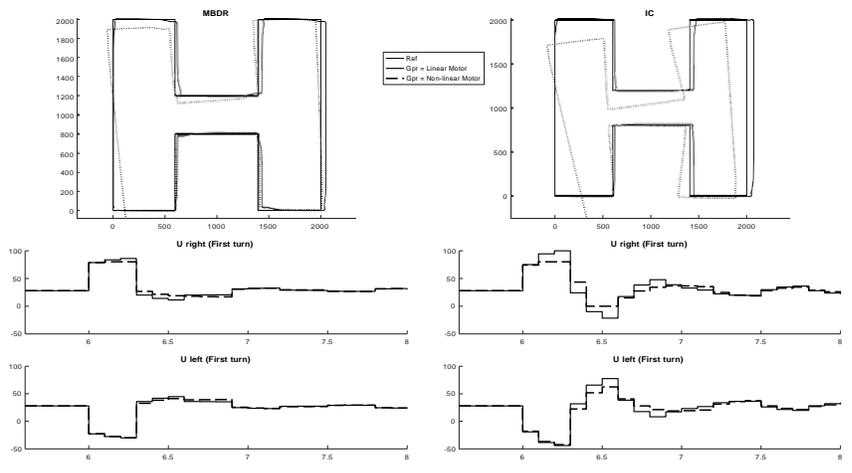

Figure 24: UGV. Linear and Non-linear servomotors comparison

41